\documentclass[10pt]{article}
\title{On the outerplanar crossing numbers of complete multipartite graphs}
\author{Adrian Riskin}
\author{\medskip Adrian Riskin\\
Department of Mathematics\\
Mary Baldwin College\\
Staunton, Virginia  24401\\
ariskin@mbc.edu}

\usepackage{amsmath,amsthm}
\theoremstyle{plain}
\newtheorem{Theorem}{Theorem}
\newtheorem{Lemma}{Lemma}
\numberwithin{Theorem}{section}
\numberwithin{Lemma}{section}
\numberwithin{equation}{section}
\begin{document}
\maketitle
\begin{abstract}
We calculate the outerplanar crossing numbers of complete multipartite graphs which have $n$ partite sets
with $m$ vertices and one partite set with $p$ vertices, where either $p|mn$ or $mn|p$.
\end{abstract}
\section{Introductory Material}

An \textit{outerplanar drawing} of a graph $G$ is a drawing of $G$ in which the vertices are 
placed on a circle and the edges are drawn as straight lines cutting through the interior.  We require
that in such drawings, no more than two edges cross in a single point.  The \textit{outerplanar
crossing number} of a graph $G$ is the minimum number of crossings taken over all outerplanar 
drawings of $G$.  We denote the outerplanar crossing number of $G$ by $\nu_{1}(G)$.
If $D$ is an outerplanar drawing of a graph $G$, the number of crossings in $D$ is denoted
by $cr_{1}(D)$.  The outerplanar crossing number of a graph
 was defined in [Kainen 1990].  There are very few exact results known.
 In fact all of them can be found in two papers:  [Fulek \textit{et al.} 2005]
  and [Riskin 2003].  Calculating these values seems to be of interest to the VLSI community, and it is 
  interesting also to graph theorists, and  
 thus we offer some new results here. Note that [Shahrokhi \textit{et al.} 1996] contains a useful
 introduction to the outerplanar crossing number problem, as well as some interesting lower bounds.
The complete $n$-partite graph with each partite set containing $m$ vertices is
denoted by $K_{m^{(n)}}$.  The complete multipartite graph $K_{\underbrace{a,\dots,a}_{m},
\underbrace{b,\dots,b}_{n}}$ is denoted by $K\left(a^{(m)},b^{(n)}\right)$. Our main results here are:
\begin{Theorem}\label{main1}
If $p|mn$ then 
\begin{equation*}
\begin{split}
\nu_{1}(K(p^{(1)},m^{(n)}))=&\frac{1}{24}m^{4}n^{4}+\frac{1}{12}m^{2}n^{3}
-\frac{1}{12}m^{4}n^{3}-\frac{1}{4}m^{3}n^{3}+\frac{1}{2}m^{3}n^{2}+\frac{1}{24}m^{4}n^{2}
-\frac{1}{4}m^{3}n\\
&+\frac{1}{6}m^{2}n^{2}p^{2}-\frac{1}{4}p^{2}mn-\frac{3}{4}m^{2}n^{2}p
+\frac{1}{12}mnp+\frac{1}{2}m^{2}np
-\frac{1}{6}m^{3}n^{2}p\\
&+\frac{1}{6}m^{3}n^{3}p
+\frac{1}{6}mn^{2}p
\end{split}
\end{equation*}
\end{Theorem}
\noindent and
\begin{Theorem}\label{main2}
If $mn|p$ then
\begin{equation*}
\nu_{1}(K(p^{(1)},m^{(n)}))=\phi(m,n,p)
-\frac{1}{12}m^{2}n^{2}\\
+\frac{1}{12}p^{2}
\end{equation*}
where $\phi(m,n,p)$ is the expression for $\nu_{1}(K(p^{(1)},m^{(n)}))$ given in Theorem \ref{main1}.
\end{Theorem}
\noindent We will need the following fact from [Riskin 2003]:
\begin{Lemma}\label{bipartitelemma}
If $m | n$ then $\nu_{1}(K_{m,n})=\frac{1}{12}n(m-1)(2mn-3m-n)$ and this minimum value is 
attained when the $m$ vertices are distributed evenly amongst the $n$ vertices.
\end{Lemma} 
\noindent Also the following, the proof of which is a mere calculation:
\begin{Lemma}\label{edgecount}
$K\left(a^{(m)},b^{(n)}\right)$ has
$ma+nb$ vertices and the number of edges is given by:
\begin{equation*}
\frac{1}{2}ma((m-1)a+nb)+\frac{1}{2}nb((n-1)b+ma)
\end{equation*}
\end{Lemma}
\noindent And finally the following from [Fulek \textit{et al.}]:
\begin{Theorem}\label{fulek}
$\nu_{1}(K_{m^{(n)}})=\frac{1}{24}m^{2}n(n-1)(m^{2}n^{2}+2n-m^{2}n-6mn+6m)$
\end{Theorem}
\section{Results}
We will need the following obvious fact, the statement of which is practically the proof:
\begin{Lemma}\label{sumlemma}
$$\sum_{k=1}^{mn}\left\lfloor\frac{k-1}{n}
\right\rfloor=\sum_{i=0}^{m-1}\sum_{k=in+1}^{(i+1)n}i=\frac{1}{2}mn(m-1)$$
\end{Lemma}
\bigskip
\noindent And we now prove our main theorems:
\begin{proof}[Proof of Theorems \ref{main1} and \ref{main2}]
Let $D$ be an outerplanar drawing of $K(p^{(1)},m^{(n)})$ with $n \geq 1$.  Denote the
$p$ vertices of the first partite set by $v_{\ell}$, $1 \leq \ell \leq p$.  There are three kinds
of crossings in $D$: First there are crossings in the drawing of $K_{m^{(n)}}$ induced by removing 
$v_{1}, \dots, v_{p}$ from $D$.  Second, there are crossings wholly in the isomorph of
$K_{p,mn}$ induced by the set of all edges joining vertices in $\{v_{1}, \dots, v_{p}\}$ to
other vertices.  Call the number of such crossings $C_{2}$.  Finally there are crossings determined by
one edge in the induced $K_{m^{(n)}}$ and one edge in the induced $K_{p,mn}$.  Call the number
of such crossings $C_{3}$.  Thus
\begin{equation*}
cr_{1}(D) \geq \nu_{1}\left(K_{m^{(n)}}\right) + C_{2} + C_{3}
\end{equation*}
Let $u_{k\ell}$ be the vertex of $K_{m^{(n)}}$ which lies $k$ spaces counterclockwise around the
circle from $v_{\ell}$.  Let $e_{k\ell}$ be the edge of $D$ joining these two vertices.  Then the 
$k-1$ vertices of $K_{m^{(n)}}$ between $u_{1\ell}$ and $u_{(k-1)\ell}$ inclusive induce
a complete multipartite graph, which we call $L_{k\ell}$.  The same holds for the $mn-k$ vertices
of $K_{m^{(n)}}$ between $u_{(k+1)\ell}$ and $u_{(mn)\ell}$ inclusive, and we refer to 
that complete multipartite graph as $R_{k\ell}$.  Let $cr_{1}(e_{k\ell})$ denote the number of 
edges of the induced $K_{m^{(n)}}$ which cross $e_{k\ell}$.  Note that 
\begin{equation*}
C_{3} = \sum_{\ell=1}^{p}\sum_{k=1}^{mn}cr_{1}(e_{k\ell})
\end{equation*}
The edges of $K_{m^{(n)}}$ which cross $e_{k\ell}$ consist of all edges of $K_{m^{(n)}}$ which
are not in $L_{k\ell}$, not in $R_{k\ell}$, and not incident to $u_{k\ell}$.  The number of 
edges incident to $u_{k\ell}$ in $K_{m^{(n)}}$ is $m(n-1)$.  Hence we can obtain a lower bound
on $cr_{1}(e_{k\ell})$ by maximizing the number of edges in $L_{k\ell}$ and in $R_{k\ell}$.
The number of edges in a complete multipartite graph with a fixed number of vertices is 
largest when the number of partite sets is as large as possible and the vertices are as evenly distributed
among the partite sets as possible.  Let $r$ be the remainder when $k-1$ is divided by $n$.  I.e.
$r=k-1-n\left\lfloor \frac{k-1}{n}\right \rfloor$.
The number of partite sets in $L_{k\ell}$ is as large as possible and the vertices are as evenly
distributed as possible amongst them when $r$ of them have $\left\lfloor\frac{k-1}{n}\right
\rfloor+1$ vertices and $n-r$ have $\left\lfloor\frac{k-1}{n}\right\rfloor$ vertices.
In other words, when 
\begin{equation*}
L_{k\ell}\cong K\left(\left(\left \lfloor \frac{k-1}{n} \right \rfloor +1\right)^
{\left(k-1-n\left\lfloor \frac{k-1}{n}\right \rfloor\right)},
\left \lfloor \frac{k-1}{n} \right \rfloor^{\left(n-k+1+n\left\lfloor \frac{k-1}{n}\right \rfloor\right)}\right)
\end{equation*}
Note that this holds even when $1 \leq k-1 \leq n$ by interpreting $K\left(1^{(r)},0^{(n-r)}\right)$
in the natural way.  A similar argument yields the fact that the number of edges in $R_{k\ell}$
is maximized when 
\begin{equation*}\label{L}
R_{k\ell}\cong K\left(\left(\left \lfloor \frac{mn-k}{n} \right \rfloor +
1\right)^{\left(mn-k-n\left\lfloor\frac{mn-k}{n}\right\rfloor\right)},
\left \lfloor \frac{mn-k}{n} \right \rfloor^{\left(n-mn+k+n\left\lfloor\frac{mn-k}{n}
\right\rfloor\right)}\right)
\end{equation*}
Then using 
Lemma \ref{edgecount} with $L_{k\ell}$ and $R_{k\ell}$, 
we find maximum values $M_{L}(k)$ and $M_{R}(k)$ for $E(L_{k\ell})$ and 
$E(R_{k\ell})$ respectively.  Hence:
\begin{equation*}\label{R}
cr_{1}(e_{k\ell}) \geq \frac{1}{2}m^{2}n(n-1)-m(n-1)-M_{L}(k)-M_{R}(k))
\end{equation*}
and therefore
\begin{align*}
C_{3}& \geq \sum_{\ell=1}^{p} \sum_{k=1}^{mn} \left[\frac{1}{2}m^{2}n(n-1)-m(n-1)-
M_{L}(k)-M_{R}(k)\right]\\
&= p \sum_{k=1}^{mn} \left[\frac{1}{2}m^{2}n(n-1)-m(n-1)-M_{L}(k)-M_{R}(k)\right]
\end{align*}
Invoking Theorem \ref{bipartitelemma} we find that 
\begin{equation*}
C_{2}\geq\nu_{1}(K_{p,mn})
=\left\{
\begin{matrix}
\frac{1}{12}mn(p-1)(2pmn-3p-mn) & p|mn\\
\enspace & \enspace\\
\frac{1}{12}p(mn-1)(2pmn-3mn-p) & mn|p
\end{matrix}
\right.
\end{equation*}
and hence
\begin{equation}\label{recursion1}
\begin{split}
\nu_{1}(K(p^{(1)},m^{(n)}))\geq \nu_{1}&(K_{m^{(n)}})+\frac{1}{12}mn(p-1)(2pmn-3p-mn)\\
&+p \sum_{k=1}^{mn} \left[\frac{1}{2}m^{2}n(n-1)-m(n-1)-M_{L}(k)-M_{R}(k)\right]
\end{split}
\end{equation}
if $p|mn$ and 
\begin{equation}\label{recursion2}
\begin{split}
\nu_{1}(K(p^{(1)},m^{(n)}))\geq \nu_{1}&(K_{m^{(n)}})+\frac{1}{12}p(mn-1)(2pmn-3mn-p)\\
&+p \sum_{k=1}^{mn} \left[\frac{1}{2}m^{2}n(n-1)-m(n-1)-M_{L}(k)-M_{R}(k)\right]
\end{split}
\end{equation}
if $mn|p$.
Applying Lemmas \ref{edgecount} and \ref{sumlemma} to the expression
\begin{equation*}
\sum_{k=1}^{mn}(M_{L}(k)+M_{R}(k))
\end{equation*}
and substituting back into \eqref{recursion1} and \eqref{recursion2} we find, after invoking 
 Theorem \ref{fulek}, the requisite expressions.
  Furthermore this bound is actually attained when the vertices of each partite set
   are distributed evenly around the circle. Incidentally, it is interesting that if $p=m$ our proof essentially reduces 
   to a proof by induction of Theorem \ref{fulek} which is different from the method used in [Fulek \textit{et al.}].
\end{proof}

\section*{Acknowledgements}
I'd like to thank Allison Ford and Georgia Weidman for helpful conversations.
\section*{Bibliography}
\begin{enumerate}
\item Fulek, R., He, H., Sykora, O., and Vrt'o, I.  Outerplanar crossing numbers of 3-row meshes, Halin graphs, 
and complete $p$-partite graphs.  Lecture Notes in Computer Science 3381(2005) 376-379.
\item Kainen, P.C. The book thickness of a graph II.  Congressus Numerantium 71(1990) pp. 127-32.
\item Riskin, A. On the outerplanar crossing numbers of $K_{m,n}$.  Bulletin of the Institute for 
Combinatorics and its Applications 39(2003) pp. 7-15.
\item Shahrokhi, F., Sz\'ekeley, L.A., S\'ykora, O., and Vr\v to, I.  The book crossing number of a graph.
J. Graph Theory 21(1996) pp. 413-24.
\end{enumerate}

\end{document}